\documentclass[letterpaper, 10 pt]{article}  

\usepackage{graphics} 
\usepackage{epsfig} 
\usepackage{amsfonts,amsmath,amsthm}
\usepackage{amssymb}

\usepackage[usenames]{color}
\definecolor{DarkBlue}{rgb}{0.1,0.1,0.5}
\definecolor{Black}{rgb}{0.0,0.0,0.0}
\definecolor{Red}{rgb}{0.9,0.0,0.1}
\definecolor{DarkBlue2}{rgb}{0.00,0.08,0.6}
\definecolor{DarkRed2}{rgb}{0.6,0.00,0.08}
\definecolor{DarkGreen2}{rgb}{0.00,0.6,0.08}

\newtheorem{thm}{Theorem}[section]%
\newtheorem{lem}{Lemma}[section]

\newtheorem{prop}{Proposition}[section]

\newcommand{\NN}{{\mathbb N}}
\newcommand{\EE}[1]{{\text{\normalsize$\mathbb E$}}\left(#1\right)}

\newcommand{\CC}{{\mathbb C}}

\newcommand{\OO}{{\mathcal O}}

\newcommand{\LL}{{\mathcal L}}

\newcommand{\II}{{\mathbb I}}

\newcommand{\tr}[1]{\text{Tr}\left(#1\right)}
\newcommand{\trr}[1]{\text{Tr}^{2}\left(#1\right)}

\newcommand{\half}{\text{\scriptsize $\frac{1}{2}$}}

\title{On stability of continuous-time quantum filters\thanks{This work was partially supported by the "Agence Nationale de la Recherche" (ANR), Projet Blanc  CQUID number 06-3-13957 and Projet Jeunes Chercheurs EPOQ2 number ANR-09-JCJC-0070.}}

\author{ H. Amini\thanks{Mines ParisTech,  Centre Automatique et Syst\`emes, Math\'{e}matiques et Syst\`{e}mes,
        60 Bd Saint Michel, 75272 Paris cedex 06, France,
        {\tt\small hadis.amini@mines-paristech.fr}}
        \and
        M. Mirrahimi \thanks{INRIA Paris-Rocquencourt,
        Domaine de Voluceau, B.P. 105, 78153 Le Chesnay cedex, France,
        {\tt\small mazyar.mirrahimi@inria.fr}}
        \and
        P.  Rouchon\thanks{Mines ParisTech,  Centre Automatique et Syst\`emes, Math\'{e}matiques et Syst\`{e}mes,
        60 Bd Saint Michel, 75272 Paris cedex 06, France,
        {\tt\small pierre.rouchon@mines-paristech.fr}}
        }

\begin{document}
\maketitle 

\begin{abstract}
We prove  that the fidelity between the quantum state governed by a continuous time  stochastic master equation driven by a Wiener process and its  associated quantum-filter state is a sub-martingale. This result is a generalization to non-pure quantum states  where  fidelity does not coincide in general with  a simple  Frobenius inner product. This result implies the stability of such filtering process but does not necessarily ensure the asymptotic convergence of such  quantum-filters.
\end{abstract}
\section{Introduction}

\medskip

The quantum filtering theory provides a foundation of statistical inference inspired in e.g. quantum optical systems. These systems are described by  continuous-time quantum stochastic differential equations. These stochastic master equations include the measurement back-action on the quantum-state.
The quantum filtering theory has been  developed by  Davies in the 1960s~\cite{davies1969quantum,davies1976quantum} and in its modern form by Belavkin in the 1980s~\cite{belavkin1995quantum,belavkin2004towards,belavkin1992quantum}.

\medskip
To these stochastic master equations are attached so-called quantum filters providing,  from the real-time measurements,  estimations of the  quantum states. Robustness and convergence of such estimation process  has been investigated in many papers. For example, sufficient convergence conditions, related to observability issues,  are given in \cite{van-handel:proba2009} and \cite{van2006filtering}. As far as we know,  general and verifiable necessary and sufficient  convergence conditions do not exist yet. For  links between quantum filtering and observers design on cones see~\cite{bonnabel-et-al:cdc2011a}.  In this paper, we generalize a stability result  for pure states (see, e.g., \cite{diosi2006coupled}) to arbitrary mixed  quantum states. More precisely, we prove that the fidelity between the quantum  state (that could be a mixed state)   and its associated quantum-filter state is a sub-martingale: this means that in average, the estimated  state  tends to be closer to the system state. This does not imply  its asymptotic convergence for large times. To prove  such  convergence, more specific analysis depending on the precise structure of the Hamiltonian, Lindbladian and measurement operators defining the system model is  required. This paper can also be seen as an extension to continuous-time evolution of~\cite{rouchon2010fidelity}.

\medskip

This paper is organized as follows. In section~\ref{sec:second}, we introduce the non linear stochastic master equations driven by Wiener processes and providing the evolutions of the quantum state and of the quantum-filter state and  we state the main result (Theorem~\ref{thm:main-tool}). Section~\ref{sec:one} is  devoted to the proof of this result:   firstly we consider an approximation via stochastic master equations driven by Poisson processes (diffusion approximation); secondly, we prove the sub-martingale  property via a time discretization. In final section, we suggest  some possible extensions of this work.

\section{Main result}\label{sec:second}

We will consider quantum systems of finite dimensions $1 <N<\infty.$ The state space of such a system is given by the set of density matrices
\begin{equation*}
{\mathcal D}:=\{\rho\in\CC^{N\times N}|\quad\rho=\rho^\dag,\quad\tr{\rho}=1,\quad\rho\geq 0\}.
\end{equation*}

Formally the evolution of the real state $\rho\in \mathcal D$ is described by the following stochastic master equation (cf. ~\cite{belavkin1992quantum,bouten2004stochastic,van2005feedback})

\begin{equation}\label{eq:master1}
d\rho_t=-\tfrac{i}{\hbar}[H,\rho_t]\;dt+\LL(\rho_t)\, dt+\,\Lambda(\rho_t)\, \;dW_t\;,
\end{equation}
where
\begin{itemize}
\item the notation $[A,B]$ refers to $AB-BA;$
\item $H=H^\dag$ is a Hermitian operator which describes the action of external forces on the system\,;
\item  $dW_t$ is the Wiener process which is the following innovation
\begin{equation}\label{eq:observation}
dW_t=dy_t-\,\tr{(L+L^\dag)\,\rho_t}\,dt\;,
\end{equation}
where $y_t$ is a continuous semi-martingale with quadratic variation $\langle y,y\rangle_t=t$ (which is the observation process obtained from the system) and $L$ is an arbitrary matrix which determines the measurement process (typically the coupling to the probe field for quantum optic systems)\;;
\item the super-operator $\LL$ is the Lindblad operator,
\begin{equation*}
\LL(\rho):=-\half\{L^\dag L,\rho\}+L \rho L^\dag,
\end{equation*}
where the notation $\{A,B\}$ refers to $AB+BA;$
\item  the super-operator $\Lambda$ is defined by
\begin{equation*}
\Lambda(\rho):= L\rho+\rho L^\dag-\tr{(L+L^\dag)\rho}\rho.
\end{equation*}
\end{itemize}
All the developments remain valid when $H$ and $L$ are deterministic  time-varying matrices. For clarity sake, we do not recall below such possible time dependence.

\medskip

The evolution of the quantum filter of state $\widehat\rho_t\in{\mathcal D}$ is  described by the following stochastic master equation which depends on the time-continuous measurement $y_t$ depending on the  true quantum state~$\rho_t$ via~\eqref{eq:observation} (see, e.g.,~\cite{barchielli1990direct}):

\begin{align}\label{eq:master2}
d\widehat{\rho}_t=-\tfrac{i}{\hbar}[H,\widehat{\rho}_t]dt&+\LL(\widehat{\rho}_t)\, dt
+\Lambda(\widehat{\rho}_t)\big(dy_t-\,\tr{(L+L^\dag)\widehat{\rho_t}}dt\big).
\end{align}

Replacing $dy_t$ by its value given in~\eqref{eq:observation}, we obtain
\begin{align*}
d\widehat{\rho}_t=&-\tfrac{i}{\hbar}[H,\widehat{\rho}_t]dt+ \LL(\widehat{\rho}_t)dt+\,\Lambda(\widehat{\rho}_t) \,dW_t\nonumber\\
&+\Lambda(\widehat{\rho}_t)\Big(\tr{(L+L^\dag)\rho_t}-\tr{(L+L^\dag)\widehat{\rho}_t}\Big)dt.
\end{align*}

A usual measurement of the difference between two quantum states $\rho$ and $\sigma$, is given by the fidelity, a real number between $0$ and $1$.  More precisely, the fidelity between $\rho$ and $\sigma$ in $\mathcal D$ is given by  (see~\cite[chapter 9]{nielsen1999quantum} for more details)
\begin{equation}\label{eq:fid}
F(\rho,\sigma)=\trr{\sqrt{\sqrt{\rho}\sigma\sqrt{\rho}}}.
\end{equation}
Here $F(\rho,\sigma)=1$ means $\rho=\sigma$, and $F(\rho,\sigma)=0$ means that the support of $\rho$ and $\sigma$ are orthogonal. $F(\rho,\sigma)$ coincides with their inner product $\tr{\rho\sigma}$ when at least one of the states $\rho$ or $\sigma$ is pure (i.e., orthogonal projector of rank  one). It is well known that the stochastic master equations~\eqref{eq:master1} and~\eqref{eq:master2} leave the domain $\mathcal D$ positively invariant.
This results form the fact that, using Ito rules, we have
\begin{multline}\label{eq:rhoJump}
 \rho_{t+dt}=\tfrac{\left(\II-\tfrac{iH}{\hbar}dt-\half L^\dag L dt+Ldy_{t}\right)~\rho_t~\left(\II-\tfrac{iH}{\hbar}dt-\half L^\dag L dt+Ldy_{t}\right)^\dag}{\tr{\left(\II-\tfrac{iH}{\hbar}dt-\half L^\dag L dt+Ldy_{t}\right)~\rho_t~\left(\II-\tfrac{iH}{\hbar}dt-\half L^\dag L dt+Ldy_{t}\right)^\dag}}
\end{multline}
and
\begin{multline}\label{eq:rhohatJump}
\widehat{\rho}_{t+dt}=\tfrac{\left(\II-\tfrac{iH}{\hbar}dt-\half L^\dag L dt+Ldy_{t}\right)~\widehat{\rho}_t~\left(\II-\tfrac{iH}{\hbar}dt-\half L^\dag L dt+Ldy_{t}\right)^\dag}{\tr{\left(\II-\tfrac{iH}{\hbar}dt-\half L^\dag L dt+Ldy_{t}\right)~\widehat{\rho}_t~\left(\II-\tfrac{iH}{\hbar}dt-\half L^\dag L dt+Ldy_{t}\right)^\dag}}
\end{multline}
where $dy_{t}= \tr{(L+L^\dag)\,\rho_{t}}\,dt + dW_{t}$.

These alternative formulations imply then directly  that, as soon as, $\rho_0$ and $\widehat\rho_0$  belong to $\mathcal D$, $\rho_t$ and $\widehat\rho_t$ remain in $\mathcal D$ for all $t \geq 0$. Therefore the expression of fidelity given by~\eqref{eq:fid} is well defined.

\medskip

\noindent We are now in position to state the main result of this paper.

\vspace{.3cm}

\begin{thm}\label{thm:main-tool} \rm  Consider the Markov processes $(\rho_t,\widehat{\rho}_t)$ satisfying the stochastic master Equations~\eqref{eq:master1} and~\eqref{eq:master2} respectively with $\rho_0$, $\widehat\rho_0$ in  $\mathcal D$. Then the fidelity $F(\rho_t,\widehat{\rho}_t)$, defined in Equation~\eqref{eq:fid}, is a submartingale, i.e.
$\EE{F(\rho_t,\widehat{\rho}_t)|(\rho_s,\widehat{\rho}_s)}\geq F(\rho_s,\widehat{\rho}_s),$ for all $t\geq s.$

\end{thm}

\medskip

We recall that the above  theorem generalize the results of~\cite{diosi2006coupled} to arbitrary purity of the real states and quantum filter.  If $\rho_0$ is pure, then $\rho_t$ remains pure for all $t>0$. In this case, $F(\rho_t,\widehat\rho_t)$ coincides with $\tr{\rho_t\widehat{\rho}_t}$. It is proved in ~\cite{diosi2006coupled} that  this Frobenius inner product is a sub-martingale for any  initial value of $\widehat{\rho}_t$:
 $\frac{d}{dt}\,\EE{\tr{\rho_t\widehat{\rho}_t}}\geq 0.$ The main idea of the proof in~\cite{diosi2006coupled} consists in using It\^{o}'s formula to reduce the theorem to showing that $\EE{\tr{d\rho_t\widehat{\rho}_t+\rho_t d\widehat{\rho}_t+d\rho_t d\widehat{\rho}_t}}\geq 0$, and then using the shift invariance of the operator $L$ in the dynamics~\eqref{eq:master1} and~\eqref{eq:master2} and choosing an appropriate value.

In the absence of any information on the purity of the real states and the quantum filter, the fidelity is given by \eqref{eq:fid}, and the application of It\^{o}'s formula for the above expression becomes much more involved. In particular, the calculation of the cross derivatives was so complicated that it became hopeless to proceed this way. As the proof presented in the next section shows, we had to choose an undirect way to approach the theorem which allowed us to avoid the heavy calculations based on second order derivative of $F$.

\section{Proof of Theorem~\ref{thm:main-tool}}\label{sec:one}

 We proceed in two steps.
\begin{itemize}
\item In the first step, we describe briefly how we obtain the stochastic master equations~\eqref{eq:master1} and~\eqref{eq:master2} as the limits of the stochastic master equations with Poisson processes using the diffusive limits inspired from the physical homodyne detection model~\cite{barchielli2009quantum,wiseman1993interpretation}.
\item In the second step, we show that the fidelity between the real state and the quantum filter which are the  solutions of  stochastic master equations with Poisson processes is a submartingale.
\end{itemize}

\medskip

\bf Step $1$. \rm  Take $\alpha>0$ a large real  number and consider the evolution of the quantum state $\rho^\alpha_t$  described by the following stochastic master equation  derived from  homodyne detection scheme (see section $6.4$ of~\cite{breuer2002theory} or~\cite{barchielli2009quantum},~\cite{wiseman1993interpretation}) for more physical details):
\begin{align}\label{eq:master3}
d\rho_t^{\alpha}=&-\tfrac{i}{\hbar}[H,\rho_t^{\alpha}]dt-\tfrac{1}{4}\Lambda_{\alpha}(\rho_t^{\alpha}) dt+\Upsilon_{\alpha}(\rho_t^{\alpha}) dN_1\\\nonumber
&-\tfrac{1}{4}\Lambda_{-\alpha}(\rho_t^{\alpha}) dt+\Upsilon_{-\alpha}(\rho_t^{\alpha}) dN_2\,,
\end{align}
where the super-operators $\Upsilon_\alpha$ is defined as follows

\begin{equation*}
\Upsilon_\alpha(\rho):= \frac{(L+\alpha)\rho (L^\dag+\alpha )}{\tr{(L+\alpha)\rho (L^\dag+\alpha )}}-\rho,
\end{equation*}

and the super-operator $\Lambda_\alpha$ is defined by
\begin{equation*}
\Lambda_\alpha(\rho):=(L+\alpha)\rho+\rho (L^\dag+\alpha )-\tr{(L+L^\dag+2\alpha) \rho}\rho.
\end{equation*}

 The super-operators $\Lambda_{-\alpha}$ and $\Upsilon_{-\alpha}$ are just obtained with replacing $\alpha$ by $-\alpha$ in the expressions given in above.

The two processes $dN_1$ and $dN_2$  are  defined by
$$dN_1:=N_{t+dt}^1-N_t^1 \quad \textrm{and} \quad dN_2:=N_{t+dt}^2-N_t^2 $$
where $N^1$ and $N^2$ are two Poisson processes. \noindent $dN_1$ and $dN_2$ take value $1$  by probabilities {\small$\half\tr{(L^\dag+\alpha )(L+\alpha)\rho_t^\alpha}dt$} and {\small$\half\tr{(L^\dag-\alpha )(L-\alpha)\rho_t^\alpha}dt$}, respectively, and take  value $0$ by the complementary probabilities.

\medskip

Similarly, the following stochastic master equation describes the infinitesimal evolution of associated   quantum filter of state $\widehat\rho^\alpha_t$ (see~\cite{barchielli1990direct}):
\begin{align}\label{eq:master4}
d\widehat{\rho_t}^{\alpha}=&-\tfrac{i}{\hbar}[H,\widehat{\rho_t}^{\alpha}]dt-\tfrac{1}{4}\Lambda_\alpha(\widehat{\rho_t}^{\alpha}) dt+\Upsilon_\alpha(\widehat{\rho_t}^{\alpha}) dN_1\\\nonumber
&-\tfrac{1}{4}\Lambda_{-\alpha}(\widehat{\rho_t}^{\alpha}) dt+\Upsilon_{-\alpha}(\widehat{\rho_t}^{\alpha}) dN_2.
\end{align}

 The following diffusive limit is obtained by the central limit theorem  when $\alpha$ tends to $+\infty$ for the semi-martingale processes applied to $dN_q$, $q=1,2$, (see ~\cite{lane1984central} or~\cite{jacod1987limit} for more details)
\begin{equation}\label{eq:limit}
dN_q\stackrel{\text{\tiny law}}{\longrightarrow} \langle\tfrac{dN_q}{dt}\rangle \,dt+\sqrt{\langle\tfrac{dN_q}{dt}\rangle}\,dW_q\,,
\end{equation}
\normalsize
where the notation $\langle A\rangle$ refers to the mean value of A. Here \\
$\langle dN_1\rangle=\half\tr{(L^\dag+\alpha )(L+\alpha)\rho_t^\alpha}\,dt$ and $\langle dN_2\rangle=\half\tr{(L^\dag-\alpha )(L-\alpha)\rho_t^\alpha}\,dt$ and  $dW_1$ and $dW_2$ are two independent Wiener processes and the convergence in~\eqref{eq:limit} is in law.

\medskip

The stochastic master Equations~\eqref{eq:master1} and~\eqref{eq:master2} are obtained by  replacing the processes $dN_q$ for $q\in\{1,2\}$ by their limits given in~\eqref{eq:limit} in the master equations~\eqref{eq:master3} and~\eqref{eq:master4} and  taking the limit when $\alpha$ goes to $+\infty$ and keeping only the lowest ordered terms in $\alpha^{-1}.$ Such a result is usually called diffusion approximation (see e.g~\cite{costantini1991diffusion}).

\medskip

Notice that {$dW$} appearing in the stochastic master equations~\eqref{eq:master1} and~\eqref{eq:master2} is given in terms of its independent constituents by

\begin{equation*}
dW=\sqrt{\half}\,\big(dW_1+dW_2\big),
\end{equation*}
and is thus itself a standard Wiener process.

 The following theorem from~\cite{pellegrini2009diffusion} justifies  the diffusion approximation described above.
 \medskip

\begin{thm}[Pellegrini-Petruccione~\cite{pellegrini2009diffusion}]\label{lem:loi}
\rm The solutions of the stochastic master Equations~\eqref{eq:master3} and~\eqref{eq:master4} converge in law, when $\alpha\rightarrow+\infty$, to the solutions of the stochastic master Equations~\eqref{eq:master1} and~\eqref{eq:master2}, respectively.
\end{thm}
\medskip

\bf Step $2$. \rm  We now prove that the fidelity between two arbitrary solutions of the stochastic master Equations~\eqref{eq:master3} and~\eqref{eq:master4} is a submartingale.
\medskip

\begin{prop}\label{prop:second}
Consider the Markov process $(\rho^{\alpha},\widehat{\rho}^{\alpha})$ which satisfy the stochastic master Equations~\eqref{eq:master3} and~\eqref{eq:master4}. Then the
fidelity defined in Equation~\eqref{eq:fid} is a
submartingale, i.e.,
\rm for all $t\geq s$, we have
\begin{equation*}
\EE{F(\rho_t^{\alpha},\widehat{\rho}_t^{\alpha})|(\rho_s^{\alpha},\widehat{\rho}_s^{\alpha})}\geq F(\rho_s^{\alpha},\widehat{\rho}_s^{\alpha}).
\end{equation*}
\end{prop}
\medskip
\begin{proof}
We consider approximations of the time-continuous Markov processes~\eqref{eq:master3} and~\eqref{eq:master4} by    discrete-time Markov processes $\xi_k$ and $\widehat{\xi_k}$:

\begin{equation}\label{eq:markov}
\xi_{k+1}=\tfrac{M_{\mu_k}\xi_k M_{\mu_k}^\dag}{\tr{M_{\mu_k}\xi_k M_{\mu_k}^\dag}}\quad\text{and}\quad \widehat{\xi}_{k+1}=\tfrac{M_{\mu_k}\widehat{\xi}_k M_{\mu_k}^\dag}{\tr{M_{\mu_k}\widehat{\xi}_k M_{\mu_k}^\dag}},
\end{equation}

where
\begin{itemize}
\item $k\in\{0,\cdots,n\}$ for a fixed large $n;$
\item initial condition $\xi_0=\rho_s^\alpha$ and $\widehat\xi_0=\widehat\rho_s^\alpha$;
\item $\mu_k$ is a random variable taking values $\mu\in\{0,1,2\}$ with probability $P_{\mu,k}=\tr{M_{\mu}\xi_k M_{\mu}^\dag};$
\item The operators $M_0,$ $M_1$ and $M_2$ are defined as follows

\begin{align*}
M_0:=1&-\tfrac{1}{4}(L^\dag+\alpha )(L+\alpha)\epsilon_n
-\tfrac{1}{4}(L^\dag-\alpha )(L-\alpha)\epsilon_n-\tfrac{i}{\hbar}H\epsilon_n;
\end{align*}
\begin{equation*}
M_1:=(L+\alpha)\sqrt{\half\epsilon_n};
\end{equation*}
and
\begin{equation*}
M_2:=(L-\alpha)\sqrt{\half\epsilon_n};
\end{equation*}
\end{itemize}
with $\epsilon_n=\tfrac{t-s}{n}.$
\medskip

In the following lemma, we show that $\xi_n$ and $\widehat{\xi}_n$ correspond to the Euler-Maruyama time discretization. Since~\eqref{eq:master3} and~\eqref{eq:master4} depend smoothly on $\rho^\alpha_t$ and $\widehat\rho^\alpha_t$, $\xi_n$ and $\widehat{\xi}_n$   converge in law towards $\rho^\alpha_t$ and $\widehat\rho^\alpha_t$ when $n\mapsto +\infty$.

\medskip

\begin{lem}\label{lem:third}
The processes $\xi_k$ and $\widehat{\xi}_k$ correspond up to second order terms in $\epsilon_n$, to the  Euler-Maruyama discretization scheme of~\eqref{eq:master3} and~\eqref{eq:master4} on  $[s,t]$.

\end{lem}

\medskip

\begin{proof} we regard the three following possible cases which arrive in according to the different values of $\mu_k$ . In each case, we show that $\xi_k$ and $\widehat{\xi}_k$ for $k\in\{0,\cdots,n\}$ are the numerical solutions of the dynamics ~\eqref{eq:master3} and~\eqref{eq:master4} respectively, with the following partition $s\leq s+\epsilon_n\leq\cdots\leq s+(n-1)\epsilon_n\leq t,$ where the uniform step length $\epsilon_n$ is $\tfrac{t-s}{n}.$
\medskip

\bf Case 1.  \rm  We first consider the case where  $\mu_k=0$ which arrives  with probability $P_{0,k}=\tr{M_0\xi_k M_0^\dag}.$ Note that

\begin{align*}
M_0\xi_k M_0^\dag=& \,\xi_k-\tfrac{1}{4}\{(L^\dag+\alpha )(L+\alpha),\xi_k\}\,\epsilon_n\\
&-\tfrac{1}{4}\{(L^\dag-\alpha )(L-\alpha),\xi_k\}\,\epsilon_n\\
&-\tfrac{i}{\hbar}[H,\xi_k]\, \epsilon_n+\OO(\epsilon_n^2).
\end{align*}
Therefore
\begin{align*}
\tr{M_0\xi_k M_0^\dag}&=1-\half\tr{(L^\dag+\alpha )(L+\alpha)\xi_k}\,\epsilon_n\\
&-\half\tr{(L^\dag-\alpha )(L-\alpha)\xi_k}\,\epsilon_n+\OO((\epsilon_n)^2)
\end{align*}

and

\begin{align*}
\big(\tr{M_0\xi_k M_0^\dag}\big)^{-1}&\approx 1+\half\tr{(L^\dag+\alpha )(L+\alpha)\xi_k}\,\epsilon_n+\\
&\half\tr{(L^\dag-\alpha )(L-\alpha)\xi_k}\,\epsilon_n+\OO((\epsilon_n)^2).
\end{align*}
Therefore, we find the following dynamics

\begin{align*}
\xi_{k+1}&\approx\xi_k-\tfrac{1}{4}\{(L^\dag+\alpha )(L+\alpha),\xi_k\}\,\epsilon_n\\
&-\tfrac{1}{4}\{(L^\dag-\alpha )(L-\alpha),\xi_k\}\,\epsilon_n\\
&+\half\tr{(L^\dag+\alpha )(L+\alpha)\xi_k}\,\xi_k\,\epsilon_n\\
&+\half\tr{(L^\dag-\alpha )(L-\alpha)\xi_k}\,\xi_k\,\epsilon_n+\OO(\epsilon_n^2).
\end{align*}
This can also be written as follows

\begin{align}\label{eq:m1}
\xi_{k+1}-\xi_k\approx&-\tfrac{1}{4}\Lambda_\alpha(\xi_k)\,\epsilon_n-\tfrac{1}{4}\Lambda_{-\alpha}(\xi_k)\,\epsilon_n+\OO(\epsilon_n^2).
\end{align}

Obviously, this dynamics in the first order of $\epsilon_n$ is equivalent to the dynamics of the numerical solution of the stochastic master Equation~\eqref{eq:master3} with the partition $s\leq s+\epsilon_n\leq\cdots\leq s+(n-1)\epsilon_n\leq t,$ when
\begin{equation*}\label{eq:case1}
N_{s+(k+1)\epsilon_n}^1-N_{s+k\epsilon_n}^1=0\quad \textrm{and}\quad N_{s+(k+1)\epsilon_n}^2-N_{s+k\epsilon_n}^2=0,
\end{equation*}

 which arrives with probability

\begin{align*}
\big(1-\half\tr{(L+\alpha)(L^\dag+\alpha )\,\xi_k}\,\epsilon_n\big)\cdots\\
\cdots\big(1-\half\tr{(L-\alpha)(L^\dag-\alpha )\,\xi_k}\,\epsilon_n\big).
\end{align*}
This probability, in the first order of $\epsilon_n$ is equal to $\tr{M_0\xi_k M_0^\dag}.$

\medskip

\bf Case 2. \rm The second case corresponds to $\mu_k=1$ which arrives with probability $\tr{M_1\xi_k M_1^\dag}.$ We find the following dynamics
\begin{equation*}
\xi_{k+1}=\tfrac{(L+\alpha)\xi_k (L^\dag+\alpha )}{\tr{(L+\alpha)\xi_k (L^\dag+\alpha )}}=\Upsilon[L+\alpha]\,\xi_k+\xi_k.
\end{equation*}

We observe that the numerical solution of the stochastic master Equation~\eqref{eq:master3} follows also the same dynamics when

 $$N_{s+(k+1)\epsilon_n}^1-N_{s+k\epsilon_n}^1=1\quad\textrm{and}\quad N_{s+(k+1)\epsilon_n}^2-N_{s+k\epsilon_n}^2=0,$$ which arrives with probability

\begin{align*}
\big(\half\tr{(L+\alpha)(L^\dag+\alpha )\,\xi_k}\,\epsilon_n\big)\big(1-\half\tr{(L-\alpha)(L^\dag-\alpha )\,\xi_k}\,\epsilon_n\big).
\end{align*}

This is equal to $\tr{M_1\xi_k M_1^\dag},$ in the first order of $\epsilon_n.$

\medskip

\bf Case 3. \rm Now we consider the last case  $\mu_k=2$ which arrives with probability $\tr{M_2\xi_k M_2^\dag}.$ Therefore, we have
\begin{equation*}
\xi_{k+1}=\tfrac{(L-\alpha)\xi_k (L^\dag-\alpha )}{\tr{(L-\alpha)\xi_k (L^\dag-\alpha )}}=\Upsilon_{-\alpha}(\xi_k)+\xi_k.
\end{equation*}
Which can also be written by the stochastic master equation~\eqref{eq:master3} with taking $\xi_k$ as the numerical solution and
$$N_{s+(k+1)\epsilon_n}^1-N_{s+k\epsilon_n}^1=0 \quad\textrm{and}\quad N_{s+(k+1)\epsilon_n}^2-N_{s+k\epsilon_n}^2=1,$$ which arrives with probability
\begin{align*}
\big(1-\half\tr{(L+\alpha)(L^\dag+\alpha )\,\xi_k}\,\epsilon_n\big)\big(\half\tr{(L-\alpha)(L^\dag-\alpha )\,\xi_k}\,\epsilon_N\big).
\end{align*}

Where in the first order of $\epsilon_n,$ this probability is equal to $\tr{M_2\xi_k M_2^\dag}.$

\medskip

Remark that, if we neglect the terms in the order of $\epsilon_n^2,$ The probability of $N_{s+(k+1)\epsilon_n}^1-N_{s+k\epsilon_n}^1=1$ and $N_{s+(k+1)\epsilon_n}^2-N_{s+k\epsilon_n}^2=1$ is negligible. Now it is clear that $\xi_k$ and similarly $\widehat{\xi}_k$ are respectively the numerical solutions of the stochastic master Equations~\eqref{eq:master3} and~\eqref{eq:master4} obtained by Euler-Maruyama method. As the right hand side of the stochastic master Equations~\eqref{eq:master3} and~\eqref{eq:master4} are smooth with respect to $\rho$ and $\widehat{\rho}$, we can use the result of~\cite[Theorem 1]{higham2005numerical} to conclude the  convergence in law  of $\xi_n$ and $\widehat{\xi}_n$ to $\rho^\alpha_t$ and $\widehat\rho^\alpha_t$ for large $n$.

\end{proof}

Now we notice that

\begin{equation*}
M_0^\dag M_0+M_1^\dag M_1+M_2^\dag M_2=\II+\OO(\epsilon_n^2):=A,
\end{equation*}
Take $\widetilde{M_r}:=(\sqrt{A})^{-1}M_r$ for $r=0,1,2$ which satisfy necessarily

\begin{equation}\label{eq:qnd}
\widetilde{M_0}^\dag\widetilde{M_0}+\widetilde{M_1}^\dag\widetilde{M_1}+\widetilde{M_2}^\dag\widetilde{M_2}=\II.
\end{equation}

Now we define the following Markov processes $\chi_k$ and $\widehat{\chi}_k$ by

\begin{equation}\label{eq:dyn1}
\chi_{k+1}=\tfrac{\widetilde{M_{\mu_k}}\chi_k\widetilde{M_{\mu_k}}^\dag}{\tr{\widetilde{M_{\mu_k}}\chi_k\widetilde{M_{\mu_k}}^\dag}}
\end{equation}
and
\begin{equation}\label{eq:dyn2}
\widehat{\chi}_{k+1}=\tfrac{\widetilde{M_{\mu_k}}\widehat{\chi}_k\widetilde{M_{\mu_k}}^\dag}{\tr{\widetilde{M_{\mu_k}}\widehat{\chi}_k\widetilde{M_{\mu_k}}^\dag}}\,,
\end{equation}
where
\begin{itemize}
\item $k\in\{0,\cdots,n\}$ for a fixed large $n;$
\item $\chi_0=\rho_s^\alpha$ and $\widehat\chi_0=\widehat\rho_s^\alpha$;
\item $\mu_k$ is a random variable taking values $\mu\in\{0,1,2\}$ with probability $P_{\mu,k}=\tr{\widetilde{M}_{\mu}\chi_k \widetilde{M}_{\mu}^\dag}.$
\end{itemize}

 Clearly $\chi_k$ and $\widehat{\chi}_k$ can also be seen as  the numerical solutions of the stochastic master Equations~\eqref{eq:master3} and~\eqref{eq:master4}, since $(\sqrt{A})^{-1}=\II-\OO(\epsilon_n^2),$ therefore in the first order of $\epsilon_n,$ the solutions $\xi_k$ and $\widehat{\xi}_k$  are equal to $\chi_k$ and $\widehat{\chi}_k,$ respectively. But, the advantage of using  $\chi_k$  and $\widehat{\chi}_k$ instead of $\xi_k$ and $\widehat{\xi}_k$ is that the operators $\widetilde{M_r}$ are Kraus operators since they satisfy Equality~\eqref{eq:qnd}. Thus we can apply Theorem~$1$ in~\cite{rouchon2010fidelity}, which proves that $F(\chi_k,\widehat{\chi}_k)$  is a sub-martingale.
 \medskip

\begin{thm}[\cite{rouchon2010fidelity}]\label{thm:discret}
\rm Consider the Markov chain $(\chi_k,\widehat{\chi}_k)$ satisfying~\eqref{eq:dyn1} and~\eqref{eq:dyn2}. Then $F(\chi_k,\widehat{\chi}_k)$ is a sub-martingale:
$
\EE{F(\chi_{k+1},\widehat{\chi}_{k+1})|(\chi_k,\widehat{\chi}_k)}\geq F(\chi_k,\widehat{\chi}_k).
$
\end{thm}
\medskip
Thus we have
$$
\EE{F(\chi_{n},\widehat{\chi}_{n})~|~ \chi_0,\widehat{\chi}_0 }\geq {F(\chi_0,\widehat{\chi}_0)}=F(\rho^\alpha_s,\widehat{\rho}^\alpha_s)
$$
Therefore by Lemma~\ref{lem:third}, we have necessarily
\begin{equation*}
\EE{F(\rho_t^{\alpha},\widehat{\rho}_t^\alpha)|\rho_s^{\alpha},\widehat{\rho}_s^\alpha)}\geq F(\rho_s^\alpha,\widehat{\rho}_s^{\alpha}),
\end{equation*}
for all $t\geq s,$ since we have (convergence in law)  $\rho_t^\alpha=\lim_{n\longrightarrow\infty}\chi_n,$  $\widehat{\rho}_t^\alpha=\lim_{n\longrightarrow\infty}\widehat{\chi}_n,$ $\chi_0=\rho_s^\alpha$ and $\widehat{\chi}_0=\widehat{\rho}_s^\alpha.$

\end{proof}
\medskip

We now apply Theorem~\ref{lem:loi} and we use the fact that  the function $F$ is bounded by one and continuous with respect to $\rho$ and $\widehat{\rho}$:
\begin{equation*}
\EE{F(\rho_t,\widehat{\rho}_t)|(\rho_s,\widehat{\rho}_s)}\geq F(\rho_s,\widehat{\rho}_s),
\end{equation*}

for all $t\geq s,$ which ends the proof of Theorem~\ref{thm:main-tool}.

\section{Numerical Test}
\begin{figure}
  \centerline{\includegraphics[width=0.7\textwidth]{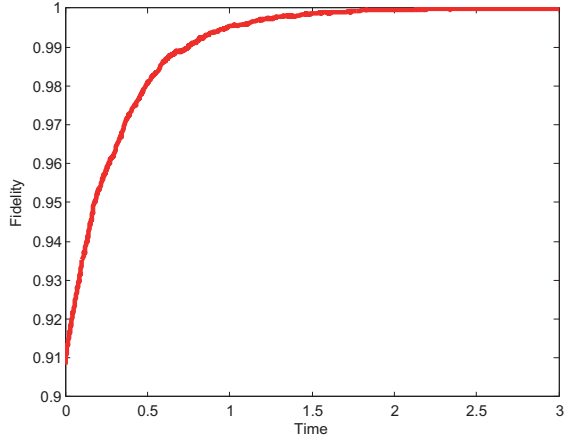}}
   \caption{The average fidelity  between the Markov processes $\rho$ and $\widehat{\rho},$ over $500$ realizations,
time $t$ from  $0$ to  $T=3$ with  discretization time step $dt=10^{-4}$.  }\label{fig:fidelity}
\end{figure}

In this section, we test the result of Theorem~\ref{thm:main-tool} through numerical simulations. Considering the two-level system of~\cite{vanHandel-et-al-2005}, we take the following Hamiltonian and measurement operators:
\begin{equation*}
H=\sigma_y=\begin{pmatrix}
0&-i \\
i&0\end{pmatrix}\qquad  L=\sigma_z=\begin{pmatrix}
1&0 \\
0&-1\end{pmatrix}.
\end{equation*}
The simulations of figure~\ref{fig:fidelity} illustrates the fidelity for 500 random trajectories starting at
\begin{equation*}
\rho_0=\begin{pmatrix}
\tfrac{1}{2}&\tfrac{1}{4} \\
\tfrac{1}{4}&\tfrac{1}{2}\end{pmatrix} \quad\textrm{and}\quad  \widehat{\rho_0}=\begin{pmatrix}
\tfrac{1}{3}&0 \\
0&\tfrac{2}{3}\end{pmatrix}.
\end{equation*}
In particular, we note that both initial states are mixed ones. As it can be seen the average fidelity is monotonically increasing. Here, the fidelity converges to one indicating the convergence of the filter towards the physical state. An interesting direction here is to characterize the situations where this convergence is ensured.

Here in order to simulate the Equations~\eqref{eq:master1} and~\eqref{eq:master2}, we have considered the alternative formulations~\eqref{eq:rhoJump} and~\eqref{eq:rhohatJump} and  the resulting  discretization scheme ($k\in\NN$ and time step  $0<dt\ll 1$)
$$
\rho_{(k+1)dt}=\tfrac{\mathcal M_k\rho_{(kdt)} \mathcal M_k^\dag}{\tr{\mathcal M_k\rho_{(kdt)} \mathcal M_k^\dag}},  \quad
\widehat{\rho}_{(k+1)dt}=\tfrac{{\mathcal M}_k\widehat{\rho}_{(kdt)}{\mathcal M}_k^\dag}{\tr{{\mathcal M}_k\widehat{\rho}_{(kdt)}{\mathcal M}_k^\dag}},
$$
where $\mathcal M_k =\II-\tfrac{iH}{\hbar}dt-\half L^\dag L dt+Ldy_{(k dt)}$ and
$ dy_{(kdt)}= \tr{(L+L^\dag)\,\rho_{(kdt)}}\,dt + dW_{(kdt)}$. For each $k$, the Wiener increment $dW_{(kdt)}$ is a centered  Gaussian random  variable of standard deviation $\sqrt{dt}$. The major interest of such  discretization is to guaranty   that, if $\rho_0,\widehat\rho_0\in\mathcal D,$ then $\rho_k$ and $\widehat\rho_k$ also remain in $\mathcal D$ for any $k\geq 0$.

\section{Concluding remarks}\label{sec:third}
The fact that the fidelity between the real quantum state and the quantum-filter state increases in average remains valid for more general stochastic master equations where other Lindbald terms are added to  $\LL(\rho)$ appearing in~\eqref{eq:master1}.  In this case the dynamics~\eqref{eq:master1} and~\eqref{eq:master2} become
\begin{multline*}
d\rho_t=-\tfrac{i}{\hbar}[H,\rho_t]\;dt+\sum_{\nu=1}^{m'}\LL'_\nu(\rho_t)\, dt
+\sum_{\mu=1}^m\LL_\mu(\rho_t)\, dt+\sum_{\mu=1}^{m}\Lambda_\mu(\rho_t)dW_t^\mu
\end{multline*}
\and
\begin{multline*}
d\widehat{\rho}_t=-\tfrac{i}{\hbar}[H,\widehat{\rho}_t]\,dt+\sum_{\nu=1}^{m'}\LL'_\nu(\widehat\rho_t)\, dt
+ \sum_{\mu=1}^m\LL_\mu(\widehat{\rho}_t)\, dt\\+\sum_{\mu=1}^m\Lambda_\mu(\widehat{\rho}_t)\bigg(dy_t^\mu-\,\tr{(L_\mu+L^\dag_\mu)\widehat{\rho_t}}dt\bigg).
\end{multline*}
where $dW_t^\mu$ are  independent Wiener processes,
\begin{equation*}
\LL_\mu(\rho):=-\half\{L_\mu^\dag L_\mu,\rho\}+L_\mu \rho L_\mu^\dag,
\end{equation*}
\begin{equation*}
\LL'_\nu(\rho):=-\half\{{L'_\nu}^\dag L'_\nu,\rho\}+L'_\nu \rho {L'_\nu}^\dag,
\end{equation*}
and $\Lambda_\mu(\rho):= L_\mu\rho+\rho L_\mu^\dag-\tr{(L_\mu+L_\mu^\dag)\rho}\rho.$

Here $m,m^\prime \geq 1$, and $(L^\prime_\nu)_{1\leq \nu\leq m^\prime}$ and $(L_\mu)_{1\leq \mu\leq m}$ are arbitrary operators. The special case  considered here corresponds to $m=1$ and $m'=1$ with $L_1=L$ and $L'_1=0.$
The formulations analogue to~\eqref{eq:rhoJump} and~\eqref{eq:rhohatJump} read then
$$
\rho_{t+dt} = \frac{(\II-dM_t) \rho_t (\II-dM_t^\dag) + \sum_{\nu=1}^{m'} L'_\nu \rho_t {L'_\nu}^\dag dt }{\tr{(\II-dM_t) \rho_t (\II-dM_t^\dag) + \sum_{\nu=1}^{m'} L'_\nu \rho_t {L'_\nu}^\dag dt}}
$$
and
$$
\widehat\rho_{t+dt} = \frac{(\II-dM_t)\widehat\rho_t (\II-dM_t^\dag) + \sum_{\nu=1}^{m'} L'_\nu \widehat\rho_t {L'_\nu}^\dag dt }{\tr{(\II-dM_t) \widehat\rho_t (\II-dM_t^\dag) + \sum_{\nu=1}^{m'} L'_\nu \widehat\rho_t {L'_\nu}^\dag dt}}
$$
where, denoting $dy_t^\mu = \tr{(L_\mu+L_\mu^\dag)\rho_t} dt + dW_t^\mu$,
$$
dM_t=\tfrac{iH}{\hbar}dt+\half\sum_{\nu=1}^{m'} {L'_\nu}^\dag L'_\nu dt 
+
\half\sum_{\mu=1}^{m} {L_\mu}^\dag L_\mu dt - \sum_{\mu=1}^{m} L_\mu dy^\mu_{t}
.
$$
For this general case, the proof of Theorem~\ref{thm:main-tool}  should follow the same lines:  first step still  relies on Theorem~\ref{lem:loi};  second  step relies now on~\cite[Theorem $2$]{rouchon2010fidelity}.

\bf Acknowledgements. \rm The authors thank C. Pellegrini, L. Zambotti and M. Gubinelli  for stimulating discussions and useful suggestions during the preparation of the paper.


\end{document}